\documentclass[preprint, 12pt]{elsarticle}

\usepackage{amsmath}
\usepackage{array}
\usepackage{fixltx2e}
\usepackage{graphicx}
\usepackage{mathtools}
\usepackage[scriptsize]{subfigure}
\usepackage{soul,color}
\DeclarePairedDelimiter{\ceil}{\lceil}{\rceil}

\newcommand{\N}{\mathit{N}}
\newcommand{\Hh}{\mathit{H}}
\newcommand{\Tna}{\mathit{T_{n,a}}}
\newcommand{\An}{\mathit{A_n}}
\journal{Special Issue on Optimization Methods in Renewable Energy Systems Design}
\begin{document}
%\linenumbers
\begin{frontmatter}
\title{Achieving an optimal trade-off between revenue and energy peak within a smart grid environment}
\author[a1]{Sezin Af\c{s}ar\corref{cor1}}
\ead{sezin.afsar@inria.fr}
\author[a1]{Luce Brotcorne}
\ead{luce.brotcorne@inria.fr}
\author[b]{Patrice Marcotte}
\ead{marcotte@iro.umontreal.ca}
\author[c]{Gilles Savard}
\ead{gilles.savard@polymtl.ca}
\address[a1]{Inria Lille-Nord Europe, Team Dolphin, 40 Avenue Halley 59650 Villeneuve d'Ascq, France}
\address[b]{DIRO, Universit\'{e} de Montr\'{e}al, C.P. 6128, Succursale Centre-Ville, Montr\'{e}al (QC) H3C 3J7, Canada}
\address[c]{D\'{e}partement de Math\'{e}matique et G\'{e}nie Industriel, \'{E}cole Polytechnique, C.P. 6079, Succursale Centre-Ville, Montr\'{e}al (QC) H3C 3A7, Canada}
\cortext[cor1]{Corresponding author. Tel: +33 59358628}

\begin{abstract}
In this paper, we consider an energy provider whose goal is to simultaneously set revenue-maximizing prices and meet a peak load constraint. The problem is cast within a bilevel setting where the provider acts as a leader (upper level) that takes into account a smart grid (lower level) that minimizes the sum of users' disutilities. The latter bases its
actions on the hourly prices set by the leader, as well as the preference schedules set by the users for each task. We consider both the monopolistic and competitive situations, and validate numerically the potential of this
approach to achieve an `optimal' trade-off between three conflicting objectives: revenue, user cost and peak demand.
\end{abstract}

\begin{keyword}
demand response, smart grid, day-ahead pricing, demand 
side management, bilevel programming, peak minimization
\end{keyword}

\end{frontmatter}

\section{Introduction}
Despite technological advancements that have allowed an increase in energy production, both traditional and green, together with a decrease of consumption in the residential and transportation markets, demand for energy is due to grow at a fast pace in the near future, putting at stress the production and distribution system, as well as the supply-demand balance. Instabilities, that trigger a chain of adverse effects for all energy users, can be mitigated by either investing to maintain a large capacity, at a high cost, or by implementing \textit{demand side management} (DSM) programs that, through controls at the customer level, make the best use of the current capacity \cite{M04,W14}. In short, DSM can be characterized by a set of tools for shaping the load curve, through peak clipping, valley filling, load shifting, strategic conservation, strategic load growth and flexible load shape \cite{G85}. In order to achieve these objectives, several programs have been put in place, such as conservation and energy efficiency programs, fuel substitution, demand response, and residential or commercial load management~\cite{M04,PD11,LMDB12}. In particular, the problem that consists of adjusting the load curve by taking explicitly into account customer reaction to prices has been addressed in several articles, such as in~\cite{MRLG10}, where residential load control through real-time
pricing has been considered. Actually, real-time pricing is frequently referred to
by economists as the most direct and efficient demand response program~\cite{B05, AE08, SMS15}. 
For a thorough literature review concerning dynamic pricing, as well as analyses of real cases, the
interested reader is referred to~\cite{BJR02}. 

In the present paper, we focus on peak load minimization through \textit{load shifting}. The importance of the issue
can be illustrated by the example of the United Kingdom, where the minimum load during summer nights is around 
30\% of the winter peak, while the average load is around 55\% of the installed generation capacity \cite{S08}, emphasizing large
fluctuations in the load curve, with close to half the generation capacity being idle for long periods of time. One foresees that the importance of load shifting will become even more apparent when the market share of plug-in cars (fully electric or hybrid) becomes significant. On average, PHEVs (plug-in electric vehicle) can be driven for 5 miles per kWh \cite{WZ11} and hence, their intensive use may double the average residential electric load \cite{MRWJSLG10}, thus putting the network at risk. 

In many countries, base load is produced by coal or nuclear power plants whereas peak load is provided by natural gas, hydro or renewable power. For this reason, electricity production during peak periods is more costly than in
the off-peak. Besides, installed production capacity has to be larger than peak load in order to ensure power supply. Since a reduction in peak load induces a decrease of production and capacity cost, it deserves to be analyzed properly. 

The specific issue that this paper addresses is Energy Peak Minimization (EPMP). It involves two decision levels: an energy provider and its customers. These two levels have conflicting objectives. The energy provider is interested in maximizing its revenue, whereas customers try to minimize their total disutility. While many articles have formulated the problem as a Nash game~\cite{MRWJSLG10,WO02,OW99,HMP00,IGCG10}, the relationship between a company and its customers better fits the leader-follower framework. More precisely, our aim is to integrate demand response explicitly into the decision making process of the energy provider. To this end, we propose a bilevel programming approach. In this setting, the \textit{leader} integrates within its decision process the reaction of the \textit{follower}. Once the leader sets his variables, the follower solves an optimization problem, taking the leader's decisions as given. Bilevel programming has been used to tackle several problems including, but 
not limited to, toll pricing~\cite{LMS98,BLMS01}, freight tariff setting~\cite{BLMS00,CLV13}, network design and pricing \cite{BCMS11,M86,BLMS08,BMS08}, electric utility planning~\cite{HN92,BJ05}. One should keep in mind that bilevel 
programs are intrinsically difficult (NP-hard) \cite{J85,HJS92}. Besides their nonconvex and combinatorial nature, the feasible region of the leader is generally nonconvex, and can be disconnected or empty~\cite{CMS05}. In the context of EPMP, bilevel programming allows to integrate DSM techniques and demand response within the optimization process of an energy provider.

The contribution of this paper is twofold. First, we develop a bilevel model for peak minimization of an energy provider, with the aim to achieve an optimal trade-off between revenue and peak power consumption without delaying demand for electricity arbitrarily. The model uses day-ahead real-time pricing and is solved for a global optimum.
We also propose a variant of the model that involves competing providers. Next, we analyze the relationship between the energy provider and its customers, where the latter are inter-connected via smart metering devices (automatic energy consumption scheduler \cite{MRWJSLG10}). In this environment, the customers not only share an energy source, but also communicate with each other via the network of smart meters which forms the smart grid. A detailed survey of smart grid and the associated enabling technologies are provided in \cite{FMXY12}. In the absence of such technology, it would be difficult for customers to keep track of hourly prices and shift their demand accordingly, 
as well as for the provider to observe the actual demand response to its pricing strategy. Smart meters enable two-way information flow and constitute an important feature of the system, allowing the application of DSM techniques. While it can be argued that metering is more expensive and may be difficult to manage for residential customers, yet it significantly decreases meter reading costs, besides assisting different pricing strategies, and improved technology provides ever easier meter management \cite{JW12}. Moreover, the smart grid system allows aggregation of residential customers with similar needs and different preferences. Since customers can act as a large client aiming at system optimum, their bargaining power is vastly increased. A survey on demand response and smart grids can be found in~\cite{S14}. 

The paper is organized as follows. The next section presents the modelling framework. It is followed by experimental results involving different parameters and instances, and a conclusion that points out 
challenges related to this field of research.
\begin{figure}[t]
\centering
\includegraphics[page=11,scale=0.33]{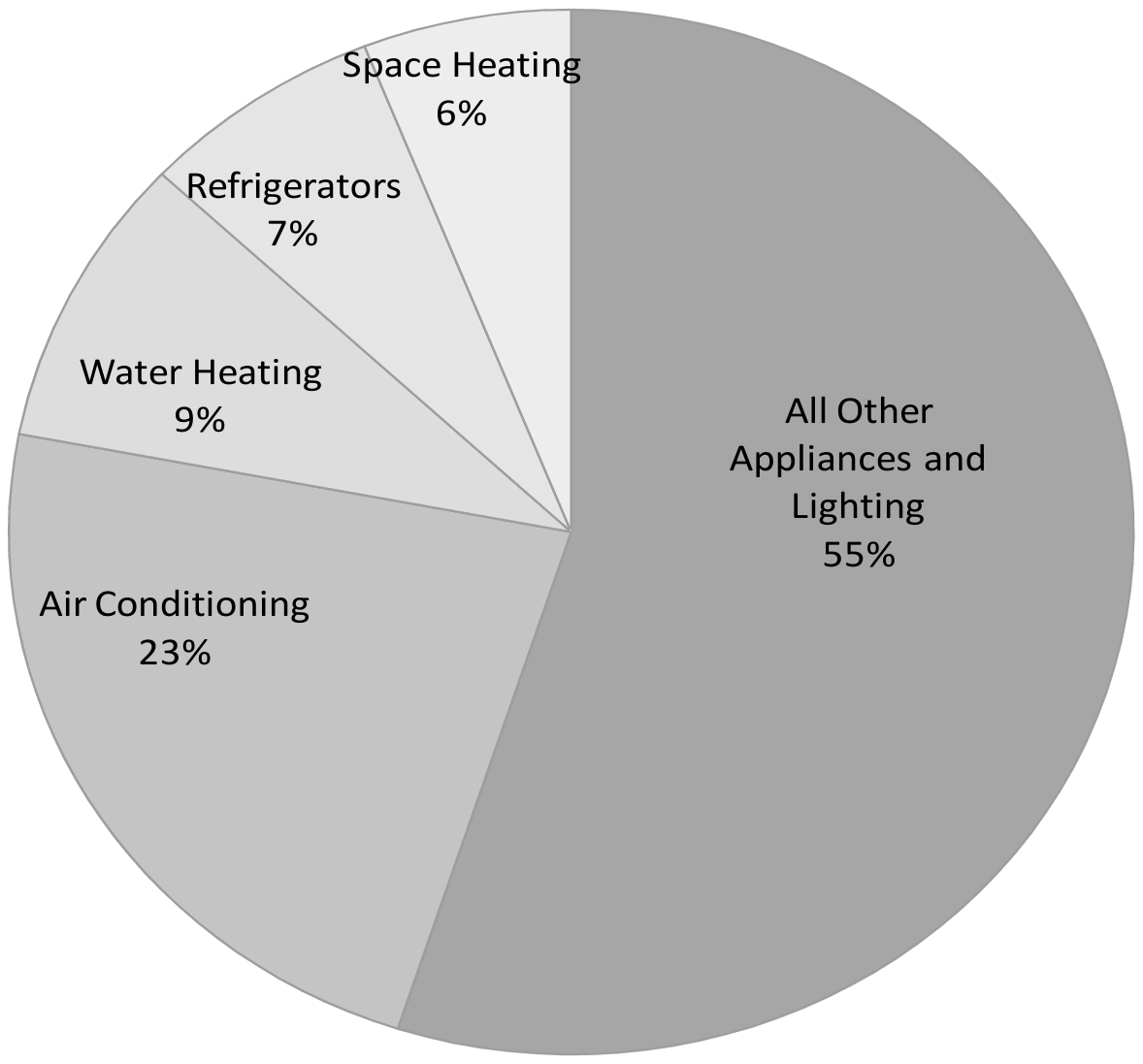}
\caption{Residential electricity use in the USA, 2011 \cite{Report12}}
\label{fig:eluse}
\end{figure}

\section{The bilevel models}
Let us consider a power sharing system involving 
a set of customers denoted by $N$, each one equipped with a smart metering device. Each customer $n$ operates a set
$A_n$ of electric residential appliances, such as air conditioners, radiators, washers, driers,
refrigerators, freezers, pool heaters, etc. The appliances can be turned on and off at any time, and their power can be adjusted at any desired level. In the United States, 45\% of household appliance consumption belongs to that category of preemptive devices (see Fig. \ref{fig:eluse}), and it follows that their intelligent control may yield a significant decrease in peak consumption.

For the sake of this study, we adopt a 24-hour planning horizon $H$. Each customer $n$ is characterized by 
its daily demands $E_{n,a}$, as well as time windows $\Tna$ transmitted to the smart meter, one per appliance
 $a\in A_n$. The set of such devices
retrieves and transmits data, and thus forms the smart grid. The grid is connected to a power source and receives hourly prices from the electricity supplier, 24 hours in advance. It allows customers to benefit from cooperation. While users are expected to have similar needs with respect to power consumption and task scheduling, some of them might yet be more reluctant to postpone their loads, whereas others are willing to switch to cheaper time slots. Such behavior differences enter the model and are taken into account by the energy provider. Precisely, the population heterogeneity with respect to price perception is captured by an \textit{inconvenience factor} specific to each customer. 

We emphasize the importance of the smart grid as a middle agent
that takes charge of scheduling, since one may not expect
the customers to monitor prices in real time, and to optimally schedule their appliances 
accordingly.

\begin{figure}[t]
\centering
\includegraphics[scale=1]{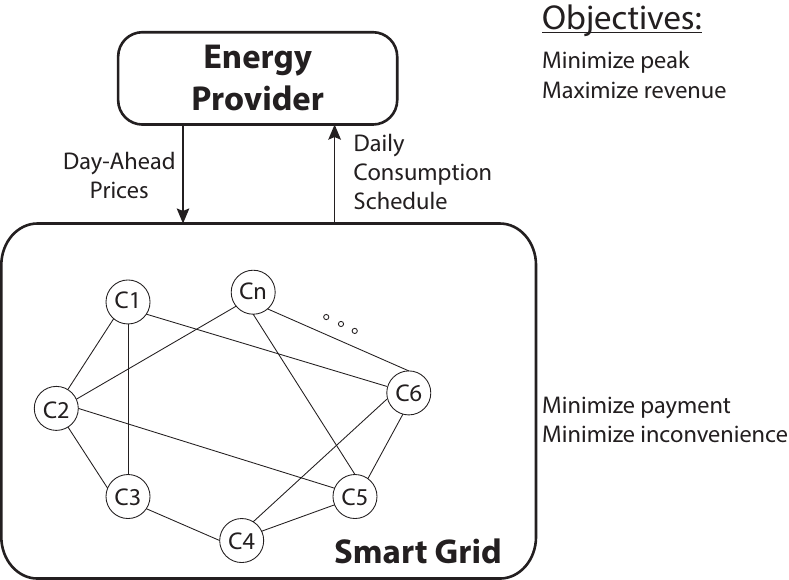}
\caption{Bilevel structure}
\label{fig:bil}
\end{figure}

In our \textit{day-ahead pricing} model, the leader applies
DSM to control and smooth out the load curve.
To tackle EPMP, two scenarios are considered: monopolistic and competitive pricing, and these will be detailed in the next two subsections. Both of them fit the bilevel paradigm, which is ideally suited at modelling game-theoretic situations
involving hierarchical players and takes the mathematical form:

\begin{equation*}
\begin{split}
\max_{x \in X} & \quad f(x,y) \\
 \text{s.t. } & \quad y\in \arg \min_{y^\prime \in Y(x)} g(x,y^\prime),
\end{split}
\end{equation*}
which usually involves two conflicting objectives. In the above model,
 the leader selects first a vector $x$, 
 taking the follower's reaction $y$ into account. In our model, the electricity supplier sets
  the prices and anticipates an optimal response from the follower. In most situations (this is the case in this paper),
 the lower level problem is convex, for given $x$, and can thus be characterized by its optimality conditions. 
 If the solution to the lower level problem is not unique, i.e., the follower is indifferent to two solutions,
 we assume that the one  most favorable to the leader is implemented. This is referred to as the `optimistic'
approach in the literature, and is justified by the fact that the solution can be made unique through
 an arbitrarily small deviation from the optimal optimistic solution. More details on this model, as well
 as the `pessimistic' (one could say `conservative') alternative, are provided in~\cite{DMZ14, DKK06}.

\subsection{Monopoly Pricing}

Let us consider a monopolistic electricity supplier, together with a group of 
customers. Once the price decision of the leader has been set, 
the smart grid automatically schedules the customer's appliances,
in compliance with the desires of the customer, i.e., it maximizes
their individual utilities.

At the upper level, the leader strives for a trade-off between revenue and peak load by maximizing
the sum of revenue minus a penalty associated with peak consumption. This is achieved by a direct control
of prices $p^h$, one for each time slot $h \in \Hh$, and indirect control of the
peak load $\Gamma$, which results from the lower level solution. 

At the lower level, the objective function involves two terms, namely \textit{electricity bill} (EB) and \textit{inconvenience cost} (IC). For a given price vector set by the leader, each individual customer minimizes weighted sum of EB and IC, which clearly conflicts with the leader's objective. 
To maximize its utility, customer $n$ selects a power level
$x_{n,a}^h$ for each appliance $a \in \An$ in each time slot $h \in \Tna$. The variable $x_{n,a}^h$ is continuous, which corresponds to preemptive jobs, and is bounded by $\beta_{n,a}^{\max}$, the power
limit of that device. Demand $E_{n,a}$ and time window $\Tna = [TW^b_{n,a},TW^e_{n,a}]$ of customer $n \in \N$ for
appliance $a \in \An$ is known and announced by the customer one day ahead. 

Throughout the paper, it is assumed that the initial time slot of a time 
window is the preferred one, and that a job cannot be performed outside its window, i.e., $x_{n,a}^h$ is defined only for $h \in \Tna$. Whenever a job is postponed within its window, a penalty that is proportional to 
the length of the delay, and inversely proportional to the width of the desired time window, is incurred. The latter assumption reflects the fact that customers that specify narrow time windows are likely to be sensitive
to the delaying of their tasks by the smart grid. More precisely, let $\lambda_n$ denote the inconvenience coefficient of a customer $n$, where a high value of $\lambda_n$ is related to a low tolerance to delay. The inconvenience cost of a job $C_{n,a}(h)$ is defined as

\begin{align*}
C_{n,a}(h) := \lambda_{n} \times E_{n,a} \times& \frac{h-TW^b_{n,a}}{TW^e_{n,a} - TW^b_{n,a}} 
&\qquad\forall n\in \N, \forall a\in \An, \forall h\in \Tna.
\end{align*}

The bilevel mathematical model of EPMP under monopolistic pricing is then expressed as:

\begin{align}
\text{(MP):}\,\,\max_{p, \Gamma} & \displaystyle \sum_{n\in \N} \sum_{a \in \An} \sum_{h\in \Tna} p^h x_{n,a}^h - \kappa \Gamma&&   \nonumber\\
\mbox{s.t. } & \displaystyle \Gamma \geq \sum_{\substack{n \in \N, a \in \An \\ \text{s.t. } h\in\Tna}} x_{n,a}^h &&\forall h \in \Hh \label{up1}\\
& \displaystyle 0 \leq p^h \leq p_{\max}^h &&\forall h \in \Hh \label{up2}\\
\min_{x} & \displaystyle \sum_{n\in \N} \sum_{a \in \An} \sum_{h\in \Tna}(p^h + C_{n,a}(h)) x_{n,a}^h &&\nonumber\\
\mbox{s.t. } & \displaystyle 0 \leq x_{n,a}^h \leq \beta_{n,a}^{\max} &&\forall {n\in \N}, \forall a\in\An, \forall h \in \Tna \label{low1} \\
& \displaystyle \sum_{h \in \Tna} x_{n,a}^h \geq E_{n,a} &&\forall {n\in \N}, \forall a\in \An, \label{low2}
\end{align}
where the `arg min' operator has been replaced by the expression of the lower level program, to simplify the exposition.

At the upper level of the above program, Constraint~(\ref{up1}) forces the variable $\Gamma$
to exceed the load in each time slot. Since it is in the interest of the leader
to minimize peak cost, $\Gamma$ should match the maximum load. 
Constraint~(\ref{up2}) sets an upper bound on the prices, and might result
from government regulation or market conditions. Constraint~(\ref{low1}) is a technical constraint that limits 
the maximum amount of power that an appliance may consume. For instance, an air conditioner can be used at most at $30^{o}$C for heating and $\beta_{n,a}^{\max}$ represents its power consumption at this level. Demand satisfaction is ensured by Constraint~(\ref{low2}).

All customers being residential users, it is natural to assume that prices only depend on time slot $h$. Note that if price discrimination were allowed, then the leader's problem would be made much easier from the computational point of view, since it would become user-separable. Unfortunately, in our model, significant modifications of the schedule may arise from a single price change, even for a single slot, thus making the problem non-trivial.

Model MP is a bilinear bilevel model involving continuous variables. Once the leader fixes his decision variables, the lower level objective function becomes linear, and is thus amenable to the classical reformulation as a single level mixed integer program (MIP) proposed by Labb\'{e} et al.~\cite{LMS98} for a network pricing problem. In this approach, the lower level's optimality conditions (primal feasibility, dual feasibility, complementary slackness) are linearized and appended to the upper level to yield an equivalent MIP formulation. In other words, the follower's mathematical program is replaced by a set of constraints that ensures the optimality of the lower level for given upper level variables. The dual and primal constraints of the follower define the feasible region of the follower, while complementary slackness constraints ensure optimality. The dual variables corresponding to constraints (\ref{low1}) and (\ref{low2}) are denoted as $w_{n,a}^h$ and $v_{n,a}$, respectively, where $w_{n,a}^h$ is defined only for $h\in \Tna$. The dual constraint corresponding to $x_{n,a}^h$ is expressed as: 

\begin{align*}
-w_{n,a}^h + v_{n,a} - p^h \leq C_{n,a}^h \qquad\forall n \in\N, \forall a\in\An, \forall h \in \Tna.
\end{align*}

Complementary slackness between $x_{n,a}^h$ and the dual constraint takes the form:
 
\begin{align*}
x_{n,a}^h(w_{n,a}^h - v_{n,a} + p^h + C_{n,a}^h)=0 \qquad \forall n \in\N, \forall a\in\An, \forall h \in \Tna
\end{align*}
and can be linearized using the fact that either $x_{n,a}^h$ or $w_{n,a}^h - v_{n,a} + p^h + C_{n,a}^h$ must be zero. Therefore, by using a sufficient large number $M_1$ and a binary variable $\psi_{n,a}^h$, one can replace the nonlinear constraint with linear ones ($\forall n \in\N, \forall a\in\An, \forall h \in \Tna$):

\begin{align*}
w_{n,a}^h - v_{n,a} + p^h + C_{n,a}^h &\leq M_1 (1-\psi_{n,a}^h) \\
x_{n,a}^h &\leq M_1 \psi_{n,a}^h  \\
\psi_{n,a}^h &\in \{0,1\} .
\end{align*}
 
Similarly, one linearizes the complementarity between dual variables and primal constraints yielding two groups of constraints. In the mixed integer linear formulation of MP, the follower's primal constraints are
given by (\ref{mip3}) and (\ref{mip6}), while the dual constraint is expressed as~(\ref{mip9}). 
Upon the introduction of binary variables $\xi_{n,a}^t$, the linearized complementary slackness between Constraint (\ref{mip3}) and $w_{n,a}^h$ corresponds to (\ref{mip4}) and (\ref{mip5}). Similarly, introducing variables $\epsilon_{n,a}$,
the complementarity between Constraint (\ref{mip6}) and $v_{n,a}$ becomes (\ref{mip7}) and (\ref{mip8}). Next,
introducing binary variable $\psi_{n,a}^t$, the complementarity between Constraints (\ref{mip9}) and $x_{n,a}^h$
is expressed as (\ref{mip10}) and (\ref{mip11}). This yields, due to the presence of identical terms (the {\it billing cost}) in the objective and the constraints, a linear expression for the leader's objective, and hence the mixed integer program

{\begin{align}
&\max_{p,\Gamma, x\atop w,v,\psi}-\sum_{\begin{subarray}{l} n \in\N \\  a\in\An\\ h \in \Tna \end{subarray}} \beta_{n,a}^{\max} w_{n,a}^h + \sum_{\begin{subarray}{l} n \in\N \\ a\in\An \end{subarray}} E_{n,a} v_{n,a}-&&\sum_{\begin{subarray}{l} n \in\N \\ a\in\An\\ h \in \Tna \end{subarray}}C_{n,a}^h x_{n,a}^h - \kappa\Gamma \label{mof}\\
\mbox{s.t.}&&& \nonumber \\
& \displaystyle \Gamma \geq \sum_{\begin{subarray}{c}n,a \\\text{s.t.} h\in \Tna \end{subarray}} x_{n,a}^h  &&\forall h \nonumber\\
 &\displaystyle 0 \leq p^h \leq p_{\max}^h  &&\forall h \in \Hh \nonumber \\
 &\displaystyle 0 \leq x_{n,a}^h \leq \beta_{n,a}^{\max}  &&\forall n \in\N,  a\in\An,  h \in \Tna  \label{mip3} \\
& \displaystyle \sum_{h\in \Tna} x_{n,a}^h \geq E_{n,a}  &&\forall n \in\N, a\in\An \label{mip6} \\
& \displaystyle -w_{n,a}^h + v_{n,a} - p^h \leq C_{n,a}^h  &&\forall n \in\N, a\in\An,  h \in \Tna \label{mip9} \\
& \displaystyle -x_{n,a}^h + M_3 \xi_{n,a}^h \leq M_3 - \beta_{n,a}^{\max} &&\forall n \in\N,  a\in\An,  h \in \Tna  \label{mip4} \\
& \displaystyle w_{n,a}^h - M_3 \xi_{n,a}^h \leq 0  && \forall n \in\N,  a\in\An,  h \in \Tna \label{mip5} \\
& \displaystyle \sum_{h\in\Tna} x_{n,a}^h + M_2 \epsilon_{n,a} \leq M_2+ E_{n,a}  &&\forall n \in\N,  a\in\An \label{mip7} \\
& \displaystyle v_{n,a} - M_2 \epsilon_{n,a} \leq 0  && \forall n \in\N,  a\in\An \label{mip8} \\
& \displaystyle w_{n,a}^h - v_{n,a} + p^h + M_1 \psi_{n,a}^h \leq M_1 - C_{n,a}^h  &&\forall n \in\N,  a\in\An, \forall h \in \Tna \label{mip10} \\
& \displaystyle x_{n,a}^h - M_1 \psi_{n,a}^h \leq 0  &&\forall n \in\N,  a\in\An,  h \in \Tna \label{mip11} \\
& \displaystyle \xi_{n,a}^h, \psi_{n,a}^h \in \{0,1\}; \,\, w_{n,a}^h \geq 0  &&\forall n \in\N, a\in\An, \forall h \in \Tna \nonumber \\
& \displaystyle  \epsilon_{n,a} \in \{0,1\}; \,\, v_{n,a} \geq 0  &&\forall n \in\N,  a\in\An. \nonumber
\end{align}}

\subsection{Competitive pricing}
Now, we turn our attention to a framework involving a competitor who declares its prices prior to the leader. The smart grid can now choose between two options for each time slot, that is, demand of each appliance can be either supplied by the leader, by the competitor or even by both of them. Therefore, the smart grid decides how much power will be supplied by the leader and the competitor in each time slot. It is important to emphasize that a customer does \textit{not} choose a supplier for all of its demand, but rather the grid chooses the \textit{total amount} of power to be purchased from each supplier.

The main concern about hour-to-hour selection of electricity supplier is payment calculation and collection. However, in a smart grid environment it is automatically calculated and billed. Therefore, several suppliers may offer prices and the most convenient one would be selected by the smart grid operator.

In the competitive pricing (CP) setting, the leader must address several issues. In addition to peak minimization and revenue maximization, as well as price ceiling constraints, market shares will only be preserved if its prices are competitive, which significantly reduces the leader's `degrees of freedom'.

In the following, we assume that competitor prices $\bar{p}^h$ are fixed, and actually assume the values $p_{\max}^h$, without loss of generality. If they were higher, then all customers would opt for the leader and the situation would be the same as in model MP. If they were lower, then $p_{\max}^h$ would be irrelevant, and leader prices would be bounded by $\bar{p}^h$. We also assume that the inconvenience factors are identical, whether electricity is supplied by the competitor or the leader. Finally, we make the conservative assumption that, whenever customers buy energy from the follower, they are automatically scheduled to their most preferred time slots.

Upon the introduction of continuous variable $\bar{x}_{n,a}^h$ to represent the amount of power purchased from the competitor, the bilevel program takes the form:

\begin{align}
\text{(CP):}\max_{p, \Gamma} & \displaystyle \sum_{n\in \N} \sum_{a \in \An} \sum_{h\in 
\Tna} p^h x_{n,a}^h - \kappa\Gamma &&\nonumber \\
\mbox{s.t. } & \displaystyle \Gamma \geq \sum_{\substack{n \in \N, a \in \An \\\text{s.t. } h\in\Tna}} x_{n,a}^h &&\forall h \nonumber \\ %\label{compup1}\\
& \displaystyle 0 \leq p^h \leq p_{\max}^h && \forall h \in \Hh \nonumber \\ %\label{compup2}\\
\min_{x,\bar{x}} &\displaystyle \sum_{n\in \N} \sum_{a \in \An} \sum_{h\in \Tna} (p^h+C_{n,a}(h))x_{n,a}^h \nonumber&&\\
+& \sum_{n\in \N} \sum_{a \in \An} \sum_{h\in \Tna} (\bar{p}^h+C_{n,a}(h))\bar{x}_{n,a}^h &&\nonumber \\
\mbox{s.t. } & \displaystyle x_{n,a}^h + \bar{x}_{n,a}^h \leq \beta_{n,a}^{\max} && \forall n \in \N ,  a \in \An,  h \in \Tna &&\label{complow1} \\
& \displaystyle \sum_{h \in \Tna}x_{n,a}^h+ \bar{x}_{n,a}^h \geq E_{n,a} && \forall n \in \N ,  \in \An \label{complow2} \\
& \displaystyle x_{n,a}^h, \bar{x}_{n,a}^h \geq 0 &&\forall n \in \N , a \in \An,  h \in \Tna, &&\label{complow3} \nonumber
\end{align}
where Constraint (\ref{complow1}) is rearranged so that the device limit is not exceeded. Constraint (\ref{complow2}) includes the amount of energy $\bar{x}_{n,a}^h$ that is purchased from the competitor, and ensures that demand is satisfied no matter which firm is the provider. 

Similar to the previous model, this program can be expressed as a single level MIP.
Following the previous notation, dual variables $w_{n,a}^h$ and 
$v_{n,a}$ are associated with (\ref{complow1}) and (\ref{complow2}), 
respectively. Next, the primal, dual and complementary slackness constraints of the lower level are
appended to the upper level, while the strong duality of the lower level is utilized 
to linearize the objective of the leader:

\begin{align*}
\max &-\sum_{\begin{subarray}{l} n \in\N \\ a\in\An\\ h \in \Tna \end{subarray}} \beta_{n,a}^{\max} w_{n,a}^h + \sum_{\begin{subarray}{l} n \in\N \\ a\in\An \end{subarray}} E_{n,a} v_{n,a} - 
\sum_{\begin{subarray}{l} n \in\N \\ a\in\An\\ h \in \Tna \end{subarray}} \bar{p}^h \bar{x}_{n,a}^h \\
&-\sum_{\begin{subarray}{l} n \in\N \\ a\in\An\\ h \in \Tna \end{subarray}} C_{n,a}(h) (x_{n,a}^h +\bar{x}_{n,a}^h) - \kappa\Gamma . 
\end{align*}

Due to the additional lower level variables $\bar{x}_{n,a}^h$, 
we must incorporate additional dual constraints, together
with their complementary slackness constraints. These are linearized as in model MP ($\forall n \in\N, \forall a \in \An,\forall h \in \Tna$):

\begin{align*}
-w_{n,a}^h + v_{n,a} &\leq C_{n,a}^h + \bar{p}^h  \\
\bar{x}_{n,a}^h \times (w_{n,a}^h - v_{n,a} + \bar{p}^h + C_{n,a}^h) &= 0 .
\end{align*} 

\section{Experimental Results and Interpretation}

In this section, we demonstrate the applicability of our approach through a number of numerical experiments involving various scenarios. The base case (BC) corresponds to setting all prices at $p_{\max}$ and
scheduling all appliances to the preferred time slots. Results of BC, MP and CP models are compared in terms of peak cost, peak load, net revenue, billing and inconvenience costs.

The scenarios involve 10 customers, each one controlling three preemptive appliances regulated
by the smart grid. The scheduling horizon is composed of 24 time slots of equal duration. 

The sensitivity of the models are tested with respect to two parameters: peak weight $\kappa$ and time window width (TWW), which is related to customer flexibility. 

Two key parameters enter the model. First, peak weight $\kappa$ reflects the importance to decrease peak load for the leader. A higher weight translates into a larger penalty, hence  a higher effort to smooth out the supply curve. Next,
time window width (TWW) provides the leader with some flexibility to induce job shifting through price adjustments, and thus indirectly smoothing out the load curve. Sensitivity of the model with respect to both parameters has been assessed, with $\kappa$ assuming values ranging from 200 to 1000, with increments of 200, and TWW
increased by 20\% or 100\% with respect to the minimum completion time 
$$\hbox{MCT} := \ceil*{{E_{n,a}}/{\beta^{\max}_{n,a}}}.$$
MCT denotes the minimum number of time slots required to meet demand $E_{n,a}$ if we could set
 all devices at their maximum level $\beta^{\max}_{n,a}$. For experimental purposes, 10 instances are randomly generated for each value of $\kappa$. In order to test TWW, similar jobs are used with different widths of time windows. 

In each scenario, parameters $\beta^{\max}_{n,a}$ and $E_{n,a}$ are uniformly generated for customer $n$ and appliance $a$. Then, the early time slot of time window for customer $n$ and appliance $a$, $TW_{n,a}^b$ is generated within 0 and $24-\ceil*{(1+TWW)\times MCT}$. The end of time window for customer $n$ and appliance $a$, $TW_{n,a}^e$ takes the value $TW_{n,a}^b+\ceil*{(1+TWW)\times MCT}$. For instance, let $\beta^{\max}$ and $E$ be equal to 2 and 8, respectively, for a given customer $n$, and appliance $a$. Then, its MCT is 4 hours. If TWW is 1.0, then its time window can start at some time slot in \{0,\ldots,16\} and must end 8 hours later. If TWW is 0.20, then $TW_{n,a}^b$ belongs to the
interval \{0,\ldots,19\} and $TW_{n,a}^e$ is $TW_{n,a}^b+5$. 

Although all customers are residential users, they may have different levels of sensitivity to delay and hence, they may behave differently. Therefore, a random inconvenience coefficient $\lambda_n$ is generated for each customer~$n$. As mentioned earlier, the inconvenience penalty function $C_{n,a}(h)$ is directly proportional with $\lambda_n$ and 
demand $E_{n,a}$, and inversely proportional to time window width. Hence, 
when $\lambda_n$ assumes a low value, customers are less delay-sensitive, 
which gives the model more flexibility to find a good schedule. Alternatively, when $\lambda_n$ assumes a large value, certain time slots become too costly and will almost never be selected. Note that in real life, $\lambda_n$ can be either selected by customers or a value can be assigned to each customer based on past data.

Both models are solved with CPLEX version 12.3 on a computer with 2.66 GHz Intel Xeon CPU and 4 GB RAM, 
running under the Windows 7 operating system. Whenever an instance could not be solved within
the time limit of 4 hours, the best integer solution has been considered. 

\begin{table}[t]
\caption{Cost Comparison of models MP and CP on 20\% TWW instances (BC = 100\%)}
\centering
{
\begin{tabular}{r r r r r r r}
$\kappa$ & MP EB & MP IC & MP TC & CP EB & CP IC & CP TC \\
\\
200 & 78.02 & 21.49 & 99.51 & 78.31 & 21.17 & 99.48 \\
400 & 77.15 & 21.59 & 98.74 & 78.01 & 20.13 & 98.14 \\
600 & 75.76 & 21.85 & 97.61 & 77.84 & 18.74 & 96.58 \\
800 & 73.52 & 22.16 & 95.68 & 78.07 & 16.81 & 94.88 \\
1000 & 71.50 & 22.48 & 93.99 & 77.63 & 16.28 & 93.91 \\
Average & 75.19 & 21.91 & 97.10 & 77.97 & 18.63 & 96.60
\end{tabular}
}
\label{tab:cost20}
\end{table}

\begin{table}[t]
\caption{Cost Comparison of models MP and CP on 100\% TWW instances (BC = 100\%)}
\centering
{
\begin{tabular}{r r r r r r r}
$\kappa$ & MP EB & MP IC & MP TC & CP EB & CP IC & CP TC \\
\\
200 & 85.12 & 14.16 & 99.28 & 85.25 & 14.06 & 99.30 \\ 
400 & 82.87 & 14.55 & 97.42 & 84.05 & 13.90 & 97.95 \\
600 & 80.13 & 15.29 & 95.42 & 84.67 & 12.69 & 97.35 \\
800 & 75.67 & 16.29 & 91.96 & 84.33 & 11.84 & 96.17 \\
1000 & 74.95 & 16.44 & 91.39 & 83.70 & 11.45 & 95.15 \\
Average & 79.75 & 15.35 & 95.10 & 84.40	& 12.79 & 97.19
\end{tabular}
}
\label{tab:cost100}
\end{table}

The first numerical results are displayed in Tables \ref{tab:cost20} and \ref{tab:cost100}.
They involve 10 random instances of 30 jobs, in both the monopolistic (MP) and competitive (MC) cases.
The user costs are split between electricity bill (EB) and inconvenience cost (IC), both percentages
being relative to the
base case (BC). For instance, the first line of Table \ref{tab:cost20} indicates that in model MP, 
out-of-pocket cost is 78.02\% and inconvenience cost is 21.49\% of the total cost 
corresponding to BC, the total (TC) being 0.5\% less than in BC, for which
the billing cost is the higher. Average values over all parameters and instances are displayed
in the final row. Models MP and CP results in a 2.9\% and 3.4\% total cost reduction, respectively for 20\% TWW instances and a 4.9\% and 2.8\% total cost reduction for 100\% instances. All values are less than 100\%, which 
reflects a cost improvement for customers for any peak weight value.

In comparison with the base case, the leader sets lower prices in order to shift some jobs to the off-peak hours, hence customers' bill is naturally reduced and inconvenience cost is increased. When peak weight $\kappa$ is large, the leader is willing to give up some revenue in order to achieve a smoother load curve. Hence, 
it lessens the bill as well. Note that EB is lower in model MP than in model CP whereas IC is higher. When the leader is a monopolist, he has to provide service to all customers. However, in the competitive case, it has the option to give up on some load in order to decrease the peak without lowering prices. According to this reasoning, IC increases as $\kappa$ increases in both tables for model MP, whereas it decreases for model CP. 

Although the total cost of the follower for 20\% TWW instances is lower in the presence of competition, it is not the case for 100\% TWW instances, and there lies an interesting fact. For instance, suppose that peak consists of a light-load and a heavy-load job alongside others, and that they are both required to be shifted in order to decrease the peak. Keeping in mind that the heavy-load job has a high inconvenience cost, the leader would have to decrease the price at least by that amount in the monopolistic case. Then, the light-load job would enjoy a price reduction that is larger than its inconvenience,
and the total cost would be lower for the light-load job and identical for the heavy one. In contrast, 
the leader can now give up on the heavy-load job in the competitive case and decrease the price only with respect to the light-load job. Hence, total cost would stay the same for both jobs. This is valid mostly for large time windows, because there are far fewer options for job shifting in 20\% instances. 

\begin{table}[t]
\caption{Comparison of models MP and CP with 20\% TWW Instances} 
\centering
{
\begin{tabular}{c r r r r r r}
 & \multicolumn{2}{c}{Av Comp Time} & \multicolumn{2}{c}{Av Gap} & 
\multicolumn{2}{c}{\# unsolved} \\ 
\multicolumn{1}{c}{($\kappa$)} & \multicolumn{1}{c}{MP} & 
\multicolumn{1}{c}{CP} & \multicolumn{1}{c}{MP} & \multicolumn{1}{c}{CP} & \multicolumn{1}{c}{MP} & \multicolumn{1}{c}{CP} \\ 
\\
200 & 1.10 & 1.20 & 0.00\% & 0.00\% & 0 & 0 \\ 
400 & 3.10 & 3.50 & 0.00\% & 0.00\% & 0 & 0 \\ 
600 & 6.80 & 13.10 & 0.00\% & 0.00\% & 0 & 0 \\ 
800 & 8.90 & 56.60 & 0.00\% & 0.00\% & 0 & 0 \\ 
1000 & 17.90 & 63.00 & 0.00\% & 0.00\% & 0 & 0 \\ 
\end{tabular}
}
\label{tab:comparison}
\end{table}

\begin{table}[t]
\caption{Comparison of models MP and CP with 100\% TWW Instances}
\centering
{
\begin{tabular}{c r r r r r r}
 & \multicolumn{2}{c}{Av Comp Time} & \multicolumn{2}{c}{Av Gap} & 
\multicolumn{2}{c}{\# unsolved} \\ 
\multicolumn{1}{c}{($\kappa$)} & \multicolumn{1}{c}{MP} & 
\multicolumn{1}{c}{CP} & \multicolumn{1}{c}{MP} & \multicolumn{1}{c}{CP} & \multicolumn{1}{c}{MP} & \multicolumn{1}{c}{CP} \\ 
\\
200 & 28.10 & 321.10 & 0.00\% & 0.00\% & 0 & 0 \\ 
400 & 339.50 & 592.67 & 0.00\% & 0.55\% & 0 & 1 \\ 
600 & 2040.20 & 1232.88 & 0.00\% & 3.87\% & 0 & 2 \\ 
800 & 4666.40 & 2350.67 & 0.00\% & 3.73\% & 0 & 4 \\ 
1000 & 7707.00 & 2034.14 & 0.00\% & 6.04\% & 0 & 3 \\ 
\end{tabular}
}
\label{tab:comparison2}
\end{table}

\begin{figure*}[h]
\centering
\subfigure[Peak Cost for 20\% TWW Instances\label{fig:peakcost20}]{\includegraphics[page=1,width=2.5in]{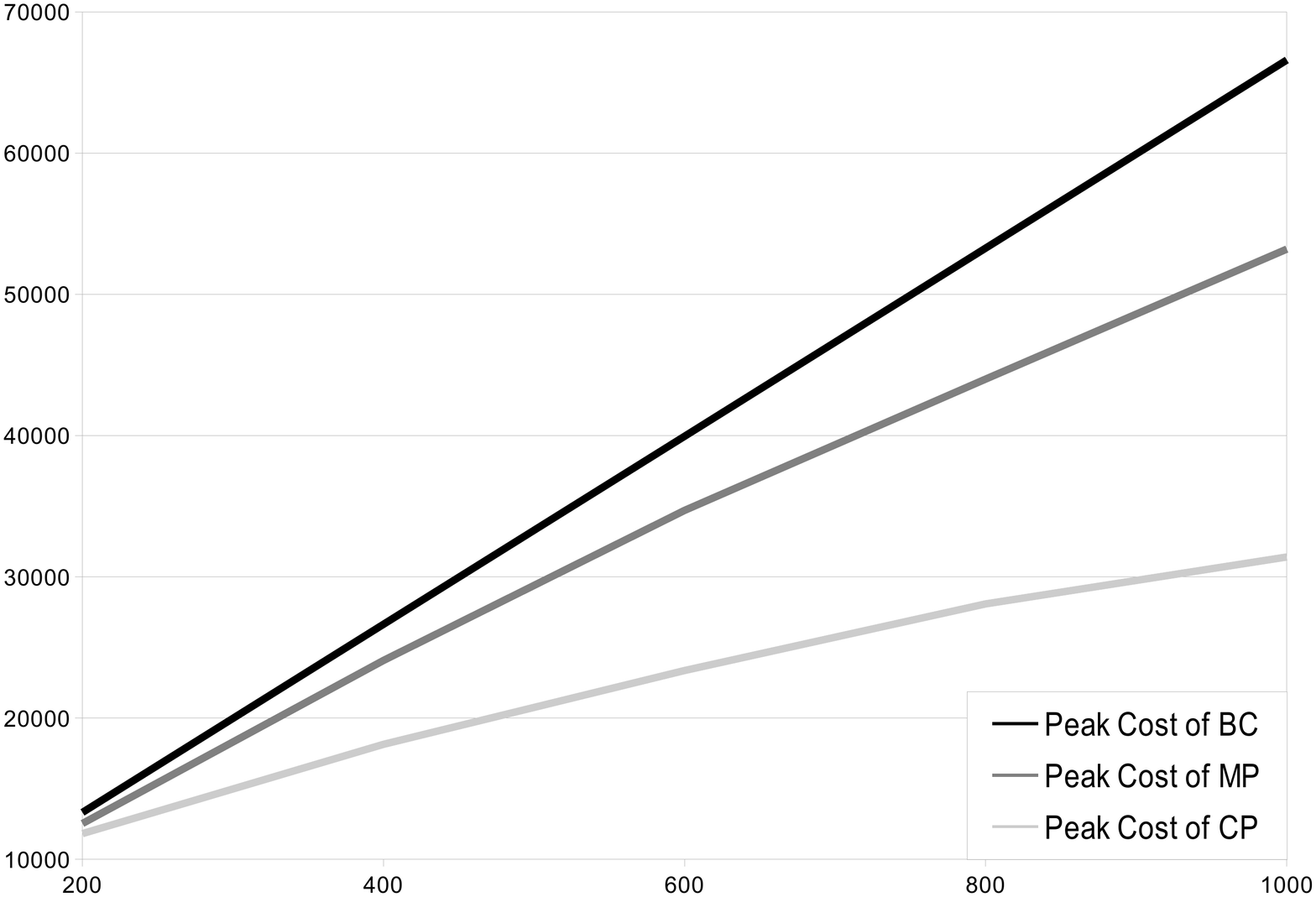}}\quad
\subfigure[Peak Cost for 100\% TWW Instances\label{fig:peakcost100}]{\includegraphics[page=2,width=2.5in]{graphics.pdf}}
\caption{Peak Cost Comparison}
\label{fig:peakcost}
\end{figure*}

\begin{figure*}[h]
\centering
\subfigure[Peak Load for 20\% TWW Instances\label{fig:peakload20}]{\includegraphics[page=3,width=2.5in]{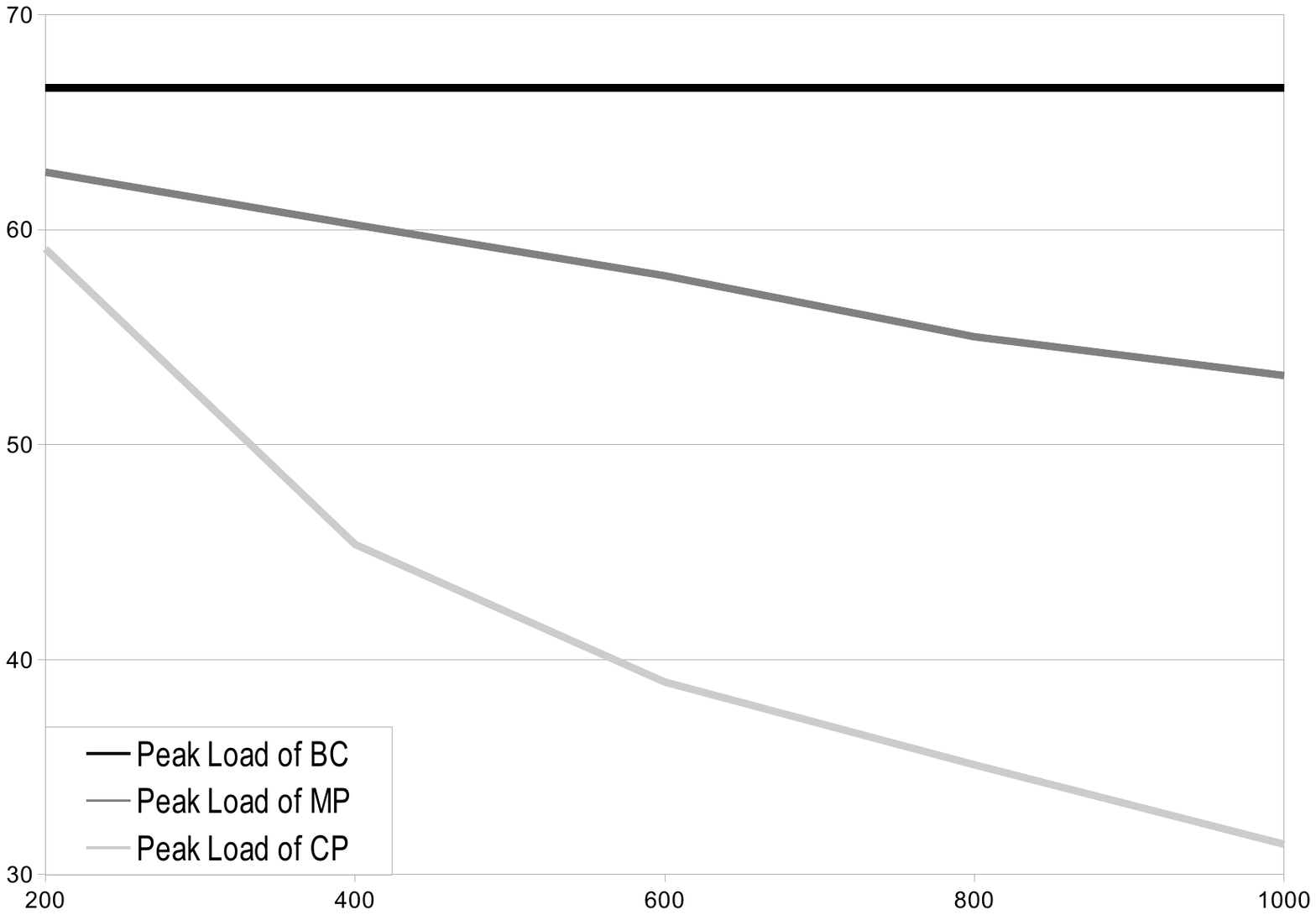}}\quad
\subfigure[Peak Load for 100\% TWW Instances\label{fig:peakload100}]{\includegraphics[page=4,width=2.5in]{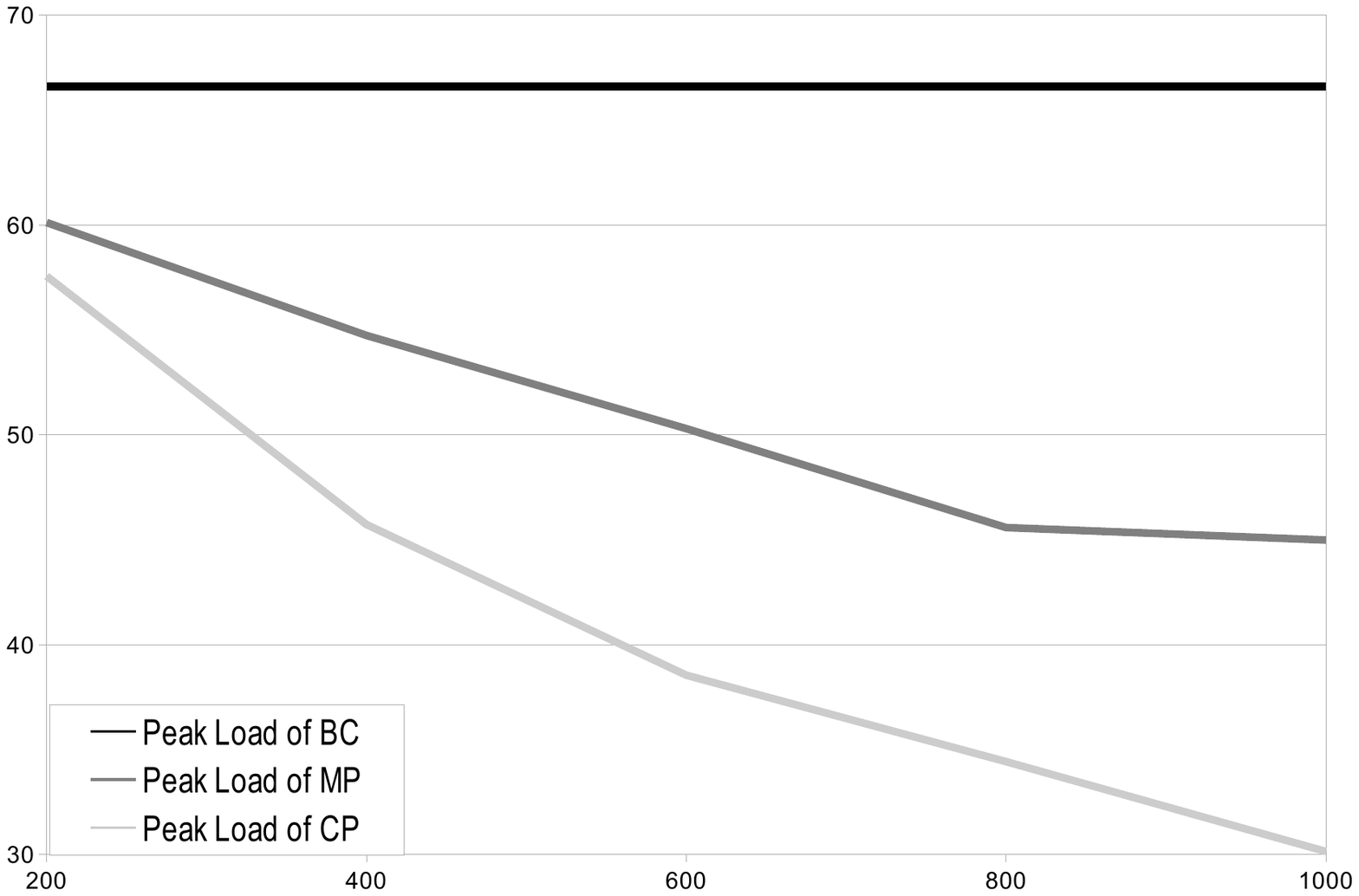}}
\caption{Peak Load Comparison}
\label{fig:peakload}
\end{figure*}

\begin{figure*}[!t]
\centering
\subfigure[Net Revenue for 20\% TWW Instances\label{fig:netrev20}]{\includegraphics[page=5,width=2.5in]{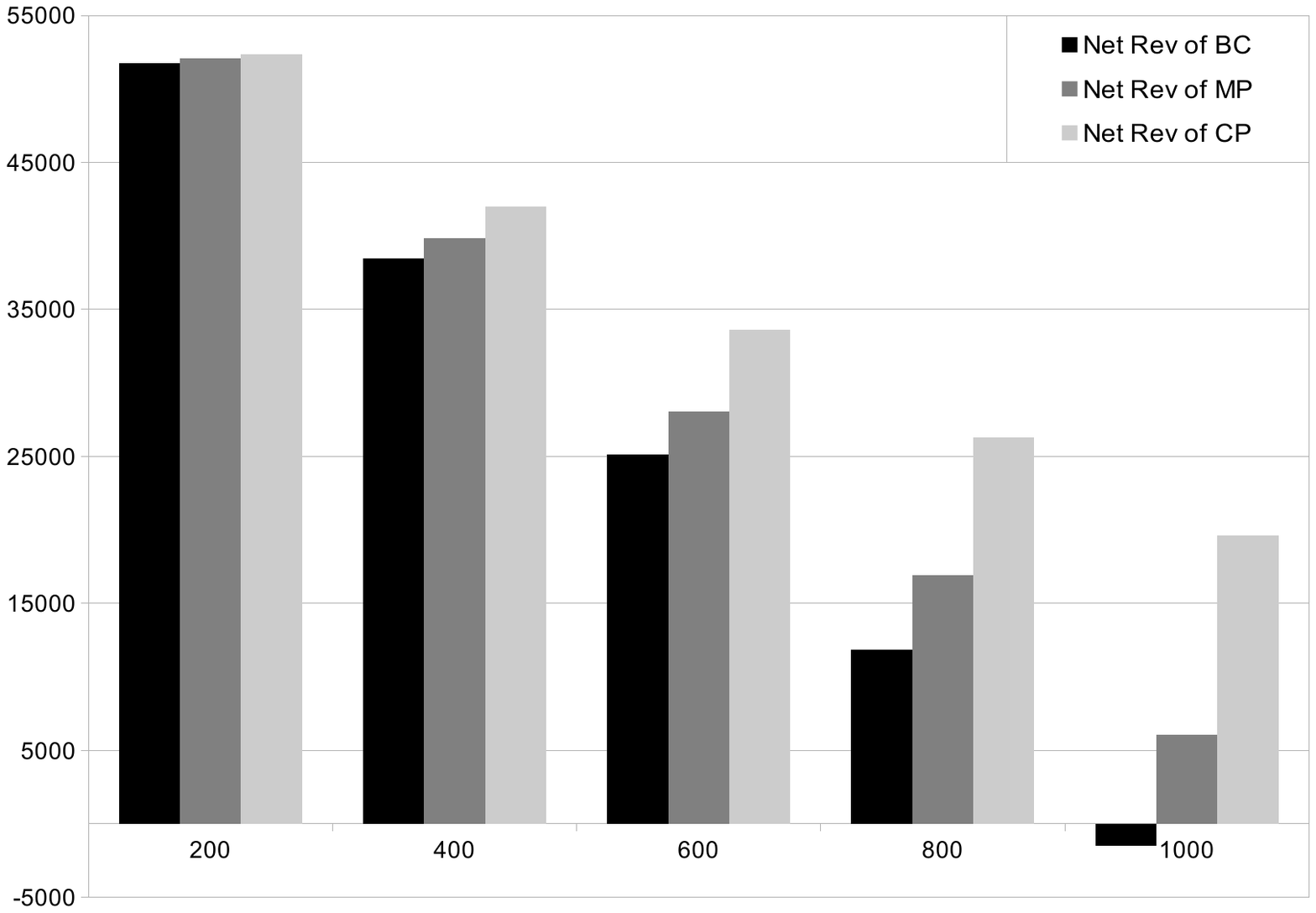}}\quad
\subfigure[Net Revenue for 100\% TWW Instances\label{fig:netrev100}]{\includegraphics[page=6,width=2.5in]{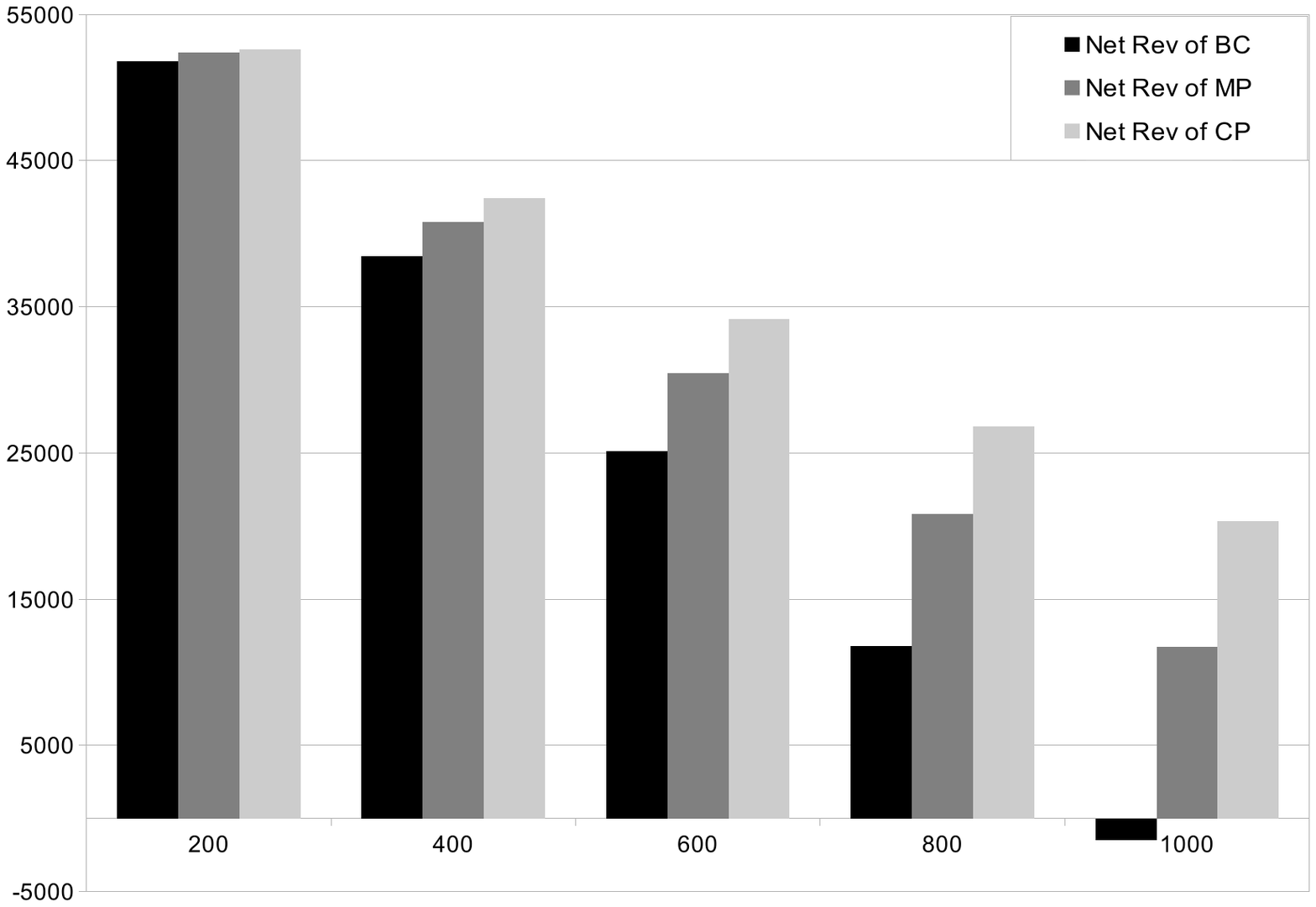}}
\caption{Net Revenue Comparison}
\label{fig:netrev}
\end{figure*}

In Tables \ref{tab:comparison} and \ref{tab:comparison2}, the computational results of model MP and CP are 
compared for 20\% and 100\% TWW instances, respectively. Average computing time (in seconds), the average 
optimality gap of unsolved instances, and the number of unsolved instances are presented. 
As aforementioned, the running time limit is set to 4 hours 
(14~400 seconds) and the values displayed in the tables are the average computation time 
over all instances that could be solved within 4 hours.

In both tables, the average computation time of both models increases together with peak weight $\kappa$,
since the leader is more willing to modify its prices in order to smooth out
the load curve, providing extra room for improvement. The average gaps and number of unsolved instances also support this argument. In addition, one can observe that average computation time is larger in Table~\ref{tab:comparison2}. Large time windows induce high running times, since there are more options to consider. Another important point is that, on average, model CP takes longer to solve if we include the unsolved instances. This result can be explained by the increased combinatorics,
the leader having the additional alternative to provide energy for a job or leave it to the competitor. In accordance with real life, competition makes the decision process more challenging.

Peak cost, peak load and net revenue (objective function value of leader) comparisons of model MP and CP to BC for 20\% instances are shown in Figures~\ref{fig:peakcost20}, \ref{fig:peakload20} and \ref{fig:netrev20}, respectively. Similar values are shown for 100\% instances in Figures~\ref{fig:peakcost100}, \ref{fig:peakload100} and \ref{fig:netrev100}. The $x$-axis consists of the peak weight parameter $\kappa$ and the 
$y$-axis represents the monetary value in Figures~\ref{fig:peakcost} and \ref{fig:netrev},
 whereas it represents peak power usage in Figure~\ref{fig:peakload}. In Figure \ref{fig:peakcost}, it can be observed that peak costs for models MP and CP increases slower than the peak cost of BC,
as weight $\kappa$ increases. Since there is a possibility of \textit{not} satisfying some of the demand for the leader in model CP, peak load and hence peak cost is lower than in model MP, as expected. In Figures \ref{fig:peakcost20} and \ref{fig:peakcost100}, it can be observed that peak cost decreases in model MP when time windows are wide, whereas it does not change much in model CP. In accordance, in Figures \ref{fig:peakload20} and \ref{fig:peakload100}, it is clear that peak load decreases in model MP. As a result, net revenue in model MP increases considerably when time windows widen. 

We now turn our attention the the leader's revenue. Average net revenue of BC is dominated by
that of model MP, and the latter is dominated by CP. Both bilevel models provide a higher \textit{net revenue} despite the
discount on some prices. In view of the peak cost constraint, the leader can adjust its
pricing strategy to increase total revenue. Perhaps more surprising, it can benefit from an open market by willingly letting demand flow to the competition, for the sake of meeting the peak constraint. It is important to note that the model behavior is very similar in both the 20\% and 100\% instances. On average, model MP provides a 13.71\% and 24.34\% net revenue increase with respect to BC on 20\% and 100\% TWW instances, respectively. Meanwhile, model CP provides a 38.31\% and 40.31\% net revenue increase with respect to BC on 20\% and 100\% TWW instances, respectively. 
 
\begin{figure*}[h]
\centering
\subfigure[Load distribution of one instance (Model MP)\label{fig:load100I}]{\includegraphics[page=7,width=2.5in]{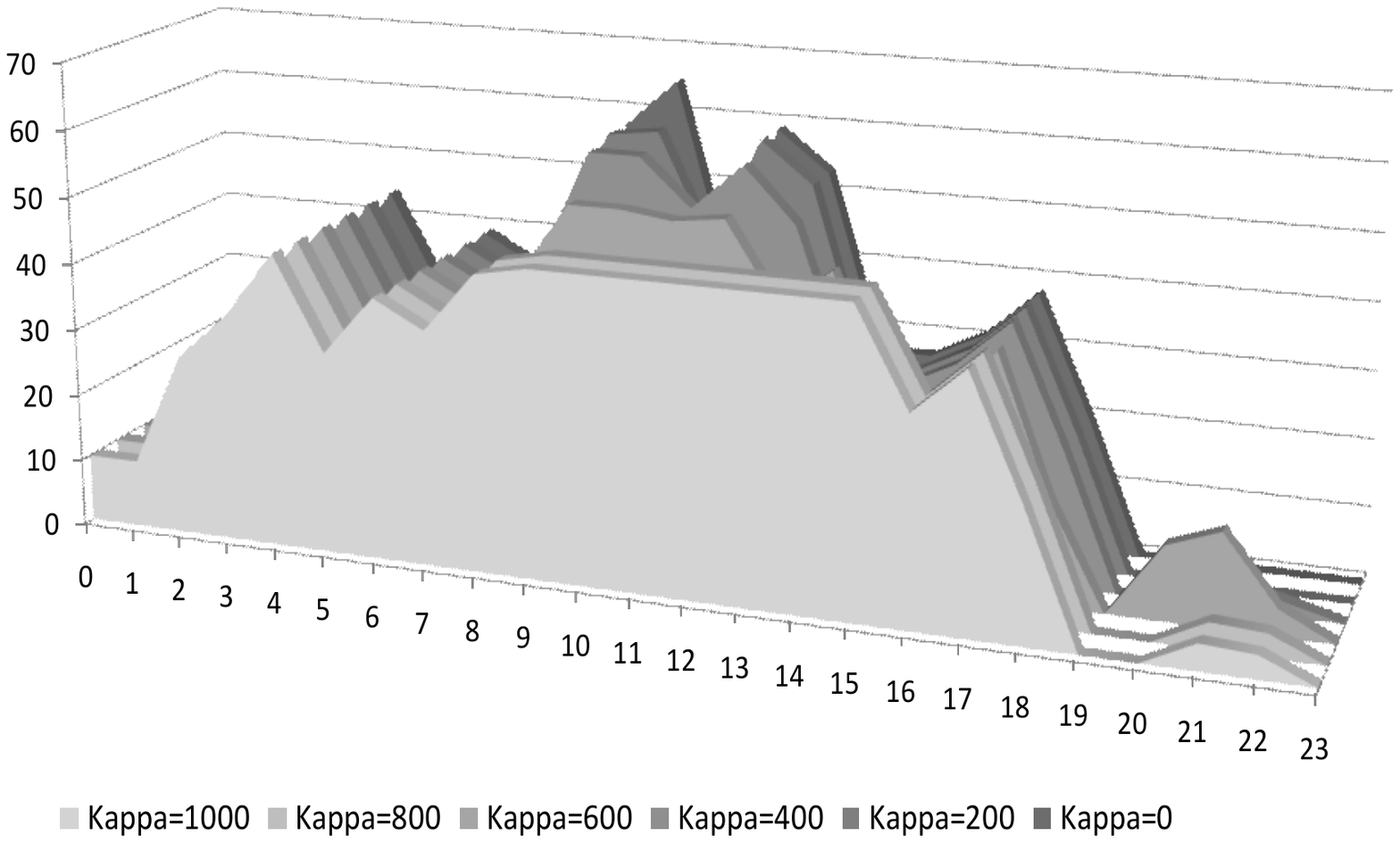}}\quad
\subfigure[Load Distribution of an Instance (Model CP)\label{fig:load100II}]{\includegraphics[page=8,width=2.5in]{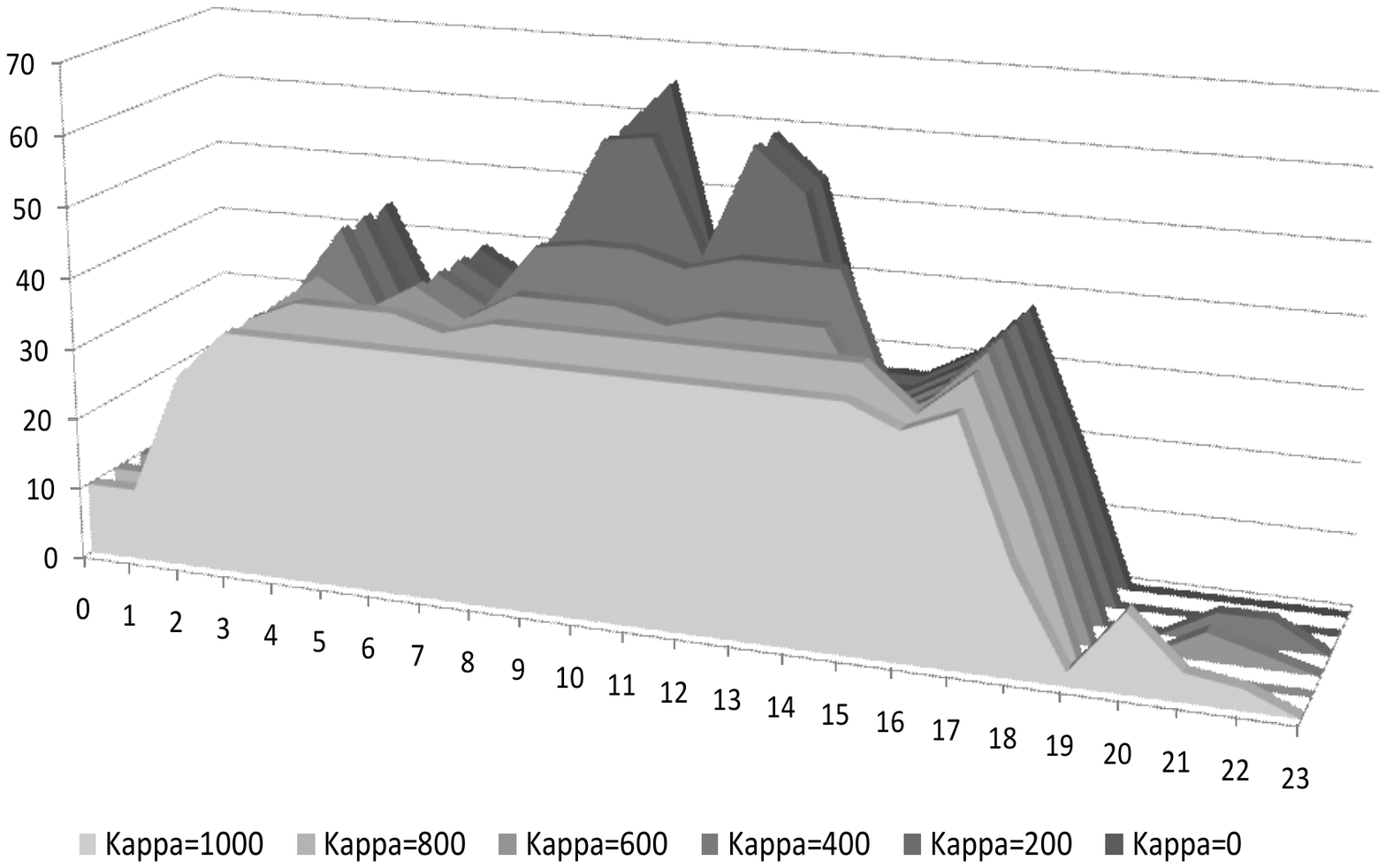}}
\caption{Load Comparison}
\label{fig:load}
\end{figure*}

\begin{figure*}[h]
\centering
\subfigure[Leader prices of one instance (Model MP)\label{fig:price100I}]{\includegraphics[page=9,width=2.5in]{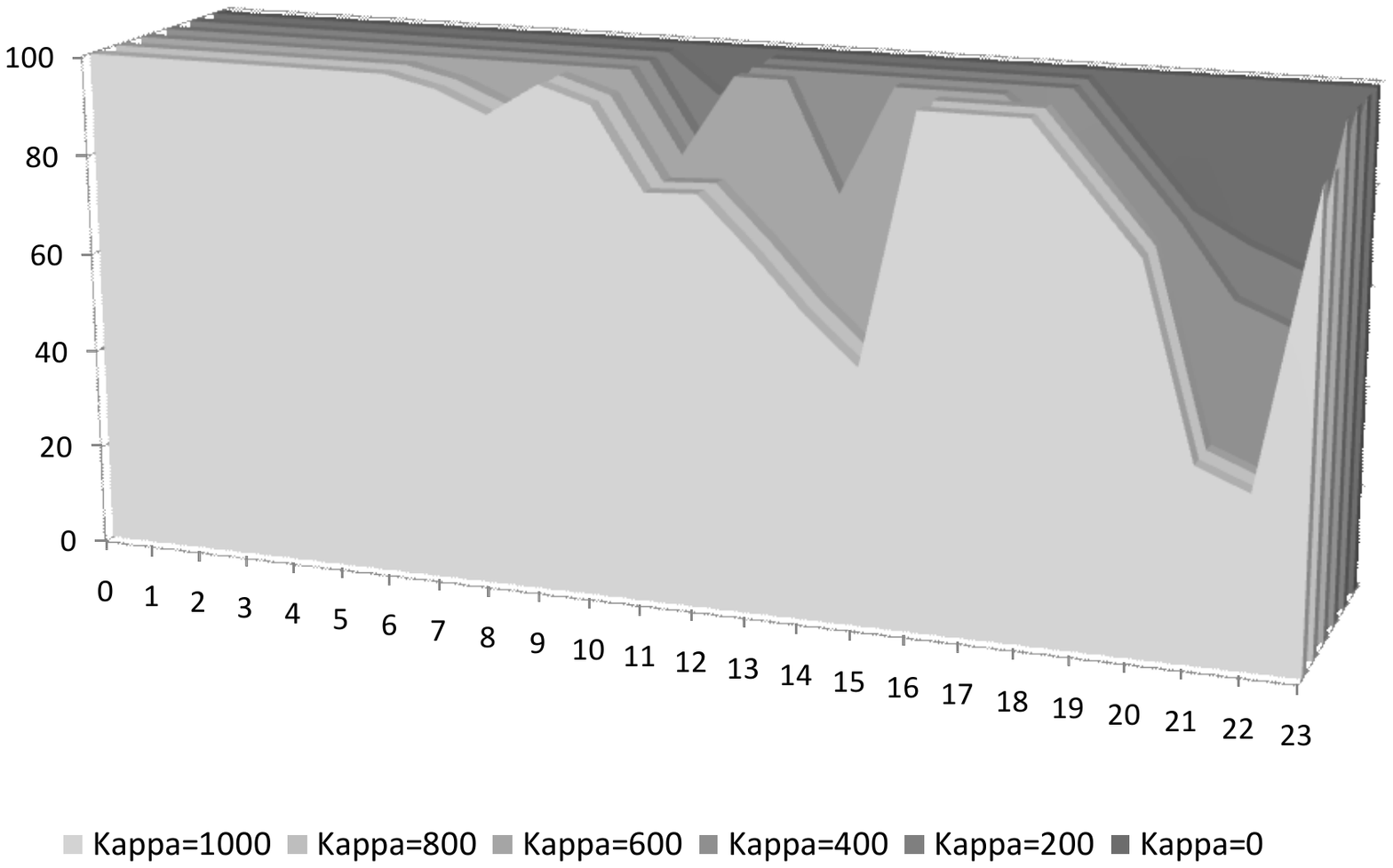}}\quad
\subfigure[Leader Prices of an Instance (Model CP)\label{fig:price100II}]{\includegraphics[page=10,width=2.5in]{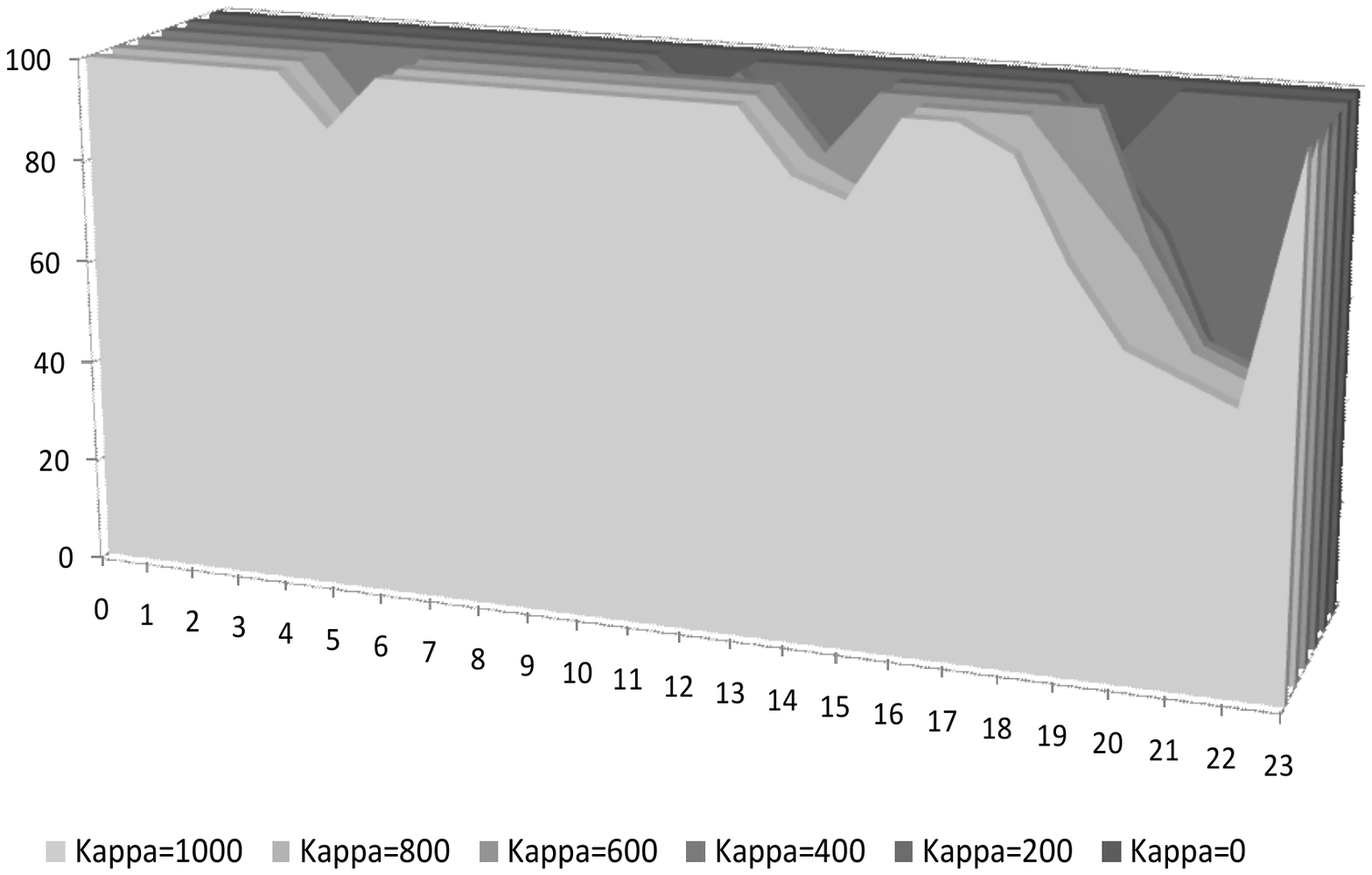}}
\caption{Price Comparison}
\label{fig:price}
\end{figure*}

Different load shapes are handled differently by the model. If the initial load curve has a single `high' peak, then it attempts to assign attractive prices around the peak to shift some of the load to later periods. However, when there are more than one peak, the shifting issue becomes more complex. Load distributions with respect to different values of $\kappa$ of an instance under model MP and 100\% TWW are shown in Figure \ref{fig:load100I}. When $\kappa$ is less than~400, the model
ignores the peak at time 12 and focuses on the one at time~10. However, as $\kappa$ increases, the model tries harder and harder to level the load curve around both peaks. In Figure~\ref{fig:price100I}, where the corresponding price vectors are displayed, low prices illustrate the effort of the leader to shift jobs around. The magnitude of price reduction escalates as $\kappa$ increases. Besides, it is again clear that prices around the other peak start moving when $\kappa$ 
exceeds the value 600. These two graphs provide a better picture of how an energy provider can achieve an optimal trade-off between revenue maximization and peak minimization.

When model CP is solved on instances of Figures \ref{fig:load100I} and \ref{fig:price100I}, load distribution and prices change as shown in Figures \ref{fig:load100II} and \ref{fig:price100II}, respectively. The leader leaves some load to the competition in return for lower peak value. By applying this strategy, it manages to keep prices higher than in the monopolistic case and achieves a smaller generation capacity. It is further observed that the leader tries to decrease peak to the level of the second highest load value (SHL). In order to achieve this, two strategies are exercised: if the time slot following peak hour has small load, then the leader opts for shifting some load to that time slot. Else, if the difference between peak load and SHL is larger than the difference between SHL and the load at the time slot following the peak, then the residual is left to the competition.

\section{Conclusion}
Maintaining an efficient supply-demand balance in the residential
energy market is a difficult task, due to variability of the demand. In order to guarantee stability during peak periods, providers must install very large capacities, or resort to costly imports. In order to address this issue, we investigated a hierarchical framework where the leader optimizes the weighted sum of revenue and peak penalty, given that the smart grid optimizes customer choice. Since electricity generation requires large investments, both for capacity building and maintenance, it is frequently handled by a private monopolist,
independent or not of the state. This corresponds to our first model. However, along with technological developments in the renewable energy generation and PHEVs, one can foresee the advent of smaller players, hence the relevance of our second model.

In this paper, we have shown that day-ahead energy pricing can be a powerful tool in terms of DSM when a smart grid system is incorporated into the system. A more efficient energy supply system is designed by flattening load curve without requiring a drastic change of habits from customers. By using pricing as a design mechanism, an optimal trade-off between
revenue and user cost can be quantified in both a monopolistic and competitive situation. By performing sensitivity with
respect to the peak load parameter, insight into efficient regulations could be achieved. 

Throughout our study, we have assumed that customers are not competing with each other,
but rather behave as a cooperative that devolves the decisions to the smart grid, providing
them with bargaining power. Analyzing EPMP under different assumptions, as well as developing 
heuristic algorithms that scale better than exact methods are two issues that will be addressed
in subsequent works.

\section*{References}
\bibliographystyle{elsarticle-num} 
 \bibliography{elsevier}
\end{document}